\documentclass[reqno,a4paper,10pt]{amsart}

\usepackage[left=2.4cm,right=2.4cm,top=3cm,bottom=3cm,a4paper]{geometry}

\usepackage{amssymb,amsfonts}
\usepackage[all,arc]{xy}
\usepackage{enumerate}

\newtheorem{thm}{Theorem}

\newtheorem{lem}{Lemma}

\theoremstyle{definition}

\theoremstyle{remark}
\newtheorem{rem}{Remark}

\title{FUNCTIONAL EQUATIONS RELATED TO THE DIRICHLET LAMBDA AND BETA FUNCTIONS}

\author{JeonWon Kim}

\begin{document}

\begin{abstract}
We give closed-form expressions for the Dirichlet beta function at even positive integers and for the Dirichlet lambda function at odd positive integers, based on the function $J(s)$ defined via convergent integral. We also show fundamental relations between Dirichlet lambda and beta functions and the function $J(s)$.
\end{abstract}

\maketitle

%\tableofcontents

\section{Introduction}  

	We will use the definitions involving the Dirichlet lambda function  and the Dirichlet beta function. 
	The Dirichlet lambda and beta function are defined as \cite{ref1}
		\begin{equation} \label{eq:1}
			\lambda(s)=\sum_{n=1}^\infty \frac{1}{(2n-1)^s}=\left(1-\frac{1}{2^s}\right)\zeta(s)   \qquad   \Re(s)>1
		\end{equation}

		\begin{equation}
			\beta(s)=\sum_{n=1}^\infty \frac{(-1)^{n-1}}{(2n-1)^s}   \qquad   \Re(s)>0
		\end{equation}
	where $\zeta(s)$ is the Riemann zeta function. The values of the Dirichlet lambda function at even positive integers and Dirichlet beta function at odd positive integers are given as \cite{ref1}

		\begin{equation}
			\lambda(2m)=(2^{2m}-1) \frac{(-1)^{m-1} \pi^{2m}}{2(2m)!}B_{2m}   \qquad   m\in \mathbb{N}
		\end{equation}

		\begin{equation}
			\beta(2m-1)=\frac{(-1)^{m-1} E_{2m-2}}{2(2m-2)!} \left(\frac{\pi}{2}\right)^{2m-1}   \qquad   m\in \mathbb{N}
		\end{equation}
	where $B_{2m}$ is Bernoulli number and $E_{2m}$ is Euler number.\\
$\newline$
	In this paper, We define the integral function $J(s)$ which can be written for all  $\Re(s)>0$
		\begin{equation}
			J(s)=\frac{1}{\Gamma(s+1)} \frac{2}{\pi} \int_0^{\frac{\pi}{2}}\frac{x^s}{\sin(x)}\,dx
		\end{equation}
	where $\Gamma$  denotes the Gamma function.
$\newline$$\newline$
	The function $J(s)$ gives closed-form expressions for the Dirichlet lambda function at odd positive integers and for the Dirichlet beta function at even positive integers.

%THEOREM 1
		
		\begin{thm} \label{thm:1}
	The values of the Dirichlet lambda function at odd positive integers are denoted by $J(s)$  as follows:
		\begin{equation}
			\lambda(2m+1)=\sum_{k=1}^m \left[ (-1)^{k-1} \lambda(2m-2k+2) \right] + (-1)^m \beta(1) J(2m)
		\end{equation} 
	where for all $m \in \mathbb{N}$
		\end{thm}

%THEOREM 2

		\begin{thm} \label{thm:2}
	The values of the Dirichlet beta function at even positive integers are denoted by $J(s)$  as follows:
		\begin{equation}
			\beta(2m)=\sum_{k=1}^n (-1)^{k-1} \beta(2m-2k+1) J(2k-1)
		\end{equation} 
	where for all $m \in \mathbb{N}$ 
		\end{thm}

$\newline$
	For example,
		\begin{gather*}
			\beta(2)=\beta(1) J(1) \\
			\beta(4)=\beta(3) J(1)-\beta(1) J(3) \\
			\beta(6)=\beta(5) J(1)-\beta(3) J(3)+\beta(1)J(5) \\
			\beta(8)=\beta(7) J(1)-\beta(5) J(3)+\beta(3)J(5)-\beta(1)J(7) \\
			\vdots 
		\end{gather*}
and
		\begin{gather*}
			\lambda(3)=\lambda(2) J(1)-\beta(1) J(2) \\
			\lambda(5)=\lambda(4) J(1)-\lambda(2) J(3)+\beta(1) J(4) \\
			\lambda(7)=\lambda(6) J(1)-\lambda(4) J(3)+\lambda(2)J(5)-\beta(1) J(6) \\
			\lambda(9)=\lambda(8) J(1)-\lambda(6) J(3)+\lambda(4)J(5)-\lambda(2)J(7)+\beta(1) J(8) \\
			\vdots 
		\end{gather*}

\section{Preliminary Lemmas}

	In this section, we start with several Lemmas used in proving Theorem $\ref{thm:1}$ and $\ref{thm:2}$.

%LEMMA 1

		\begin{lem} \label{lem:1}
	If $n$ is a positive integer, then  \cite{ref2}
		\begin{align}
			\sum_{k=1}^n \cos((2k-1)x)&=\frac{1}{2}\csc(x)\sin(2nx) \\
			\sum_{k=1}^n \sin((2k-1)x)&=\csc(x) \sin^2(nx)
		\end{align}
		\end{lem}
\
		\begin{proof}
	Consider the following sum,
			\begin{equation*}
				S=\sum_{k=1}^n \cos((2k-1)x)+i\sum_{k=1}^n \sin((2k-1)x)=\sum_{k=1}^n e^{i(2k-1)x}
			\end{equation*}
	Since $S$ is a geometric series with common ratio $e^{2ix}$ 
			\begin{align*}
				S&=\frac{e^{ix}(1-e^{2nix})}{1-e^{2ix}}=\frac{(e^{-nix}-e^{nix})e^{nix}}{e^{-ix}-e^{ix}}=\frac{-2i\sin(nx)\{\cos(nx)+i\sin(nx)\}}{-2i\sin(x)} \\
				&=\frac{1}{2} \csc(x)\sin(2nx)+i\csc(x)\sin^2(nx)
			\end{align*}	 
	Therefore,

		\begin{align*}
			\sum_{k=1}^n \cos((2k-1)x)&=\frac{1}{2}\csc(x)\sin(2nx) \\
			\sum_{k=1}^n \sin((2k-1)x)&=\csc(x) \sin^2(nx)
		\end{align*}
		\end{proof}

%LEMMA 2

		\begin{lem} \label{lem:2}
	If $n$ is a positive integer, then  \cite{ref2}
		\begin{align}
			\sum_{k=1}^n (-1)^{k-1}\cos((2k-1)x)=\sec(x)\sin^2\left(\frac{n(\pi-2x)}{2}\right) 
		\end{align} 
		\end{lem}
\
		\begin{proof}
	Consider the following sum,
			\begin{align*}
				S&=\sum_{k=1}^n (-1)^{k-1}\cos((2k-1)x)+i\sum_{k=1}^n (-1)^{k-1}\sin((2k-1)x)=\sum_{k=1}^n (-1)^{k-1} e^{i(2k-1)x} \\
				&=\frac{e^{ix}(1-(-1)^{n}e^{2nix})}{1+e^{2ix}}=\frac{(1-(-1)^{n}e^{2nix})}{e^{-ix}+e^{ix}}=\frac{1-(-1)^{n}\cos(2nx)-(-1)^{n}i\sin(2nx)}{2\cos(x)}
			\end{align*}
	Taking the real part,
			\begin{equation*}
				\Re(S)=\frac{1-(-1)^n\cos(2nx)}{2\cos(x)}=\frac{1-\cos(n\pi)\cos(2nx)}{2\cos(x)}=\frac{1-cos(n\pi-2nx)}{2\cos(x)}=\frac{\sin^2(n(\pi-2x)/2)}{\cos(x)}
			\end{equation*}
	Therefore,
		\begin{align*}
			\sum_{k=1}^n (-1)^{k-1}\cos((2k-1)x)=\sec(x)\sin^2\left(\frac{n(\pi-2x)}{2}\right) 
		\end{align*} 
		\end{proof}

%LEMMA 3

		\begin{lem} \label{lem:3}
	Let $A$ be a $n\times n$ matrix, and $ A_{ij}=\sin\left(\dfrac{(2i-1)(2j-1)\pi}{4n}\right)$, then $A^{-1}=\dfrac{2}{n}A$
		\end{lem}
\
		\begin{proof}
	Note that  $(i, j)$th element of the matrix $A^2$. The $A^2$ is the $n\times n$ matrix whose $(i, j)$th entry is given by
			\begin{equation*}
				A_{ij}^2=\sum_{m=1}^n \left[\sin\left(\frac{(2i-1)(2m-1)\pi}{4n}\right)
					\sin\left(\frac{(2j-1)(2m-1)\pi}{4n}\right) \right]
			\end{equation*}
	If $i=j$, we have
			\begin{equation*}
				A_{ij}^2=\sum_{m=1}^n \sin^2\left(\frac{(2i-1)(2m-1)\pi}{4n}\right)\
			\end{equation*}
	By using the identity $\sin^2(x)=\frac{1}{2}(1-\cos(2x))$ and Lemma $\ref{lem:1}$.
			\begin{align*}
				A_{ij}^2&=\sum_{m=1}^n \left[\frac{1}{2}-\frac{1}{2}\cos\left(\frac{(2i-1)(2m-1)\pi}{2n}\right)\right]
					=\frac{n}{2}-\frac{1}{2}\sin((2i-1)\pi)\csc\left(\frac{(2i-1)\pi}{2n}\right) \\
					&=\frac{n}{2}
			\end{align*}
	If $i\neq j$, we have
			\begin{align*}
				A_{ij}^2&=\frac{1}{2}\sum_{m=1}^n \left[\cos\left(\frac{(2i-2j)(2m-1)\pi}{4n}\right)
					-\cos\left(\frac{(2i+2j-2)(2m-1)\pi}{4n}\right) \right] \\
					&=\frac{1}{4}\sin((i-j)\pi)\csc\left(\frac{(j-i)\pi}{2n}\right)
					-\frac{1}{4}\sin((i+j-1)\pi)\csc\left(\frac{(i+j-1)\pi}{2n}\right) \\
					&=0
			\end{align*}
	Thus, if $i=j$, the expression evaluates to $n/2$ and if $i\neq j$, the this expression evaluates to 0. By the two cases above,
			\begin{equation*}
				A^2=\frac{n}{2} I_n
			\end{equation*}
	where $I_n$ is $n\times n$ identity matrix. Therefore $A$ is non-singular and
			\begin{equation*}
				A^{-1}=\frac{2}{n}A
			\end{equation*}
		\end{proof}

%LEMMA 4

		\begin{lem} \label{lem:4}
	Let $B$ be a $n\times n$ matrix, and $ B_{ij}=\cos\left(\dfrac{(2i-1)(2j-1)\pi}{4n}\right)$, then $B^{-1}=\dfrac{2}{n}B$
		\end{lem}
\
		\begin{proof}
	Note that  $(i, j)$th element of the matrix $B^2$. The $B^2$ is the $n\times n$ matrix whose $(i, j)$th entry is given by
			\begin{equation*}
				B_{ij}^2=\sum_{m=1}^n \left[\cos\left(\frac{(2i-1)(2m-1)\pi}{4n}\right)
					 \cos\left(\frac{(2j-1)(2m-1)\pi}{4n}\right) \right]
			\end{equation*}
	If $i=j$, we have
			\begin{equation*}
				B_{ij}^2=\sum_{m=1}^n \cos^2 \left( \frac{(2i-1)(2m-1)\pi}{4n}\right) 
			\end{equation*}
	By using the identity $\cos^2(x)=\frac{1}{2}(1+\cos(2x))$ and Lemma $\ref{lem:1}$.
			\begin{align*}
				B_{ij}^2&=\sum_{m=1}^n \left[\frac{1}{2}+\frac{1}{2}\cos\left(\frac{(2i-1)(2m-1)\pi}{2n}\right)\right]
					=\frac{n}{2}+\frac{1}{2}\sin((2i-1)\pi)\csc\left(\frac{(2i-1)\pi}{2n}\right) \\
					   &=\frac{n}{2}
			\end{align*}
	If $i\neq j$, we have
			\begin{align*}
				B_{ij}^2&=\frac{1}{2}\sum_{m=1}^n \left[\cos\left(\frac{(2i-2j)(2m-1)\pi}{4n}\right)
					+\cos\left(\frac{(2i+2j-2)(2m-1)\pi}{4n}\right)\right] \\
					&=\frac{1}{4}\sin((i-j)\pi)\csc\left(\frac{(j-i)\pi}{2n}\right)
					+\frac{1}{4}\sin((i+j-1)\pi)\csc\left(\frac{(i+j-1)\pi}{2n}\right) \\
					&=0
			\end{align*}
	Finally, the expression for $i=j$ evaluates to $n/2$, and the expression for $i\neq j$ evaluates to 0. By the two cases above,
			\begin{equation*}
				B^2=\frac{n}{2} I_n
			\end{equation*}
	where $I_n$ is $n\times n$ identity matrix. Therefore $B$ is non-singular and
			\begin{equation*}
				B^{-1}=\frac{2}{n}B
			\end{equation*}
		\end{proof}

%LEMMA 5

		\begin{lem} \label{lem:5}
	Let $f(s)$ be an infinite series defined by
			\begin{equation}
				f(s)=\dfrac{1}{\Gamma(s+1)} \lim_{n \to \infty} \dfrac{1}{n} \sum_{p=1}^n 
					\dfrac{\left( \dfrac{(2p-1)\pi}{4n} \right)^s}{\sin \left(\dfrac{(2p-1)\pi}{4n} \right)} 
			\end{equation}
	where $\Re(s)>0$, then $f(s)=J(s)$.
		\end{lem}
\
		\begin{proof}
	$f(s)$ is represented by difference of two infinite series as follows:
			\begin{equation*}
				f(s)=\frac{1}{\Gamma(s+1)} \lim_{m \to \infty} \frac{1}{m} \sum_{k=1}^{2m} 
					\frac{\left( \left(\dfrac{\pi}{2}\right) \dfrac{k}{2m} \right)^s}
						{\sin\left( \left(\dfrac{\pi}{2}\right) \dfrac{k}{2m} \right)} 
					-\frac{1}{\Gamma(s+1)} \lim_{n \to \infty} \frac{1}{n} \sum_{k=1}^{n}  
					\frac{\left( \left(\dfrac{\pi}{2}\right) \dfrac{2k}{2n} \right)^s}
						{\sin\left( \left(\dfrac{\pi}{2}\right) \dfrac{2k}{2n} \right)} 
			\end{equation*}
	By substituting $2m=n$,
			\begin{align*}
				f(s)&=\frac{1}{\Gamma(s+1)} \lim_{n \to \infty} \frac{2}{n} \sum_{k=1}^{n} 
					\frac{\left( \left(\dfrac{\pi}{2}\right) \dfrac{k}{n} \right)^s}
						{\sin\left( \left(\dfrac{\pi}{2}\right) \dfrac{k}{n} \right)}
					-\frac{1}{\Gamma(s+1)} \lim_{n \to \infty} \frac{1}{n} \sum_{k=1}^{n} 
					 \frac{\left( \left(\dfrac{\pi}{2}\right) \dfrac{2k}{2n} \right)^s} 
						 {\sin\left( \left(\dfrac{\pi}{2}\right) \dfrac{2k}{2n} \right)}  \\
					&=\frac{1}{\Gamma(s+1)} \lim_{n \to \infty} \frac{1}{n} \sum_{k=1}^{n}  
					\frac{\left( \left(\dfrac{\pi}{2}\right) \dfrac{k}{n} \right)^s}
						{\sin\left( \left(\dfrac{\pi}{2}\right) \dfrac{k}{n} \right)} 
			\end{align*}
	Let $\Delta x=\big(\frac{\pi}{2} \big) \frac{1}{n}$, $x_{k}=\big(\frac{\pi}{2} \big) \frac{k}{n}$, $f(x)=\frac{x^s}{\sin(x)}$, then

			\begin{align*}
				f(s)&=\frac{1}{\Gamma(s+1)} \lim_{n \to \infty} \frac{2}{\pi} \sum_{k=1}^n f(x_{k})
					\Delta x= \frac{1}{\Gamma(s+1)} \frac{2}{\pi} \int_{0}^{\frac{\pi}{2}} f(x) dx \\
				     &=\frac{1}{\Gamma(s+1)} \frac{2}{\pi} \int_{0}^{\frac{\pi}{2}} \frac{x^s}{\sin(x)} dx = J(s)
			\end{align*}
		\end{proof}

%LEMMA 6

		\begin{lem} \label{lem:6}
	Let $W(s)$ be a divergent function defined by
			\begin{equation}
				W(s)=\frac{1}{\Gamma(s+1)} \lim_{n \to \infty} \frac{1}{n} \sum_{p=1}^n 
					 \frac{\left( \dfrac{(2p-1)\pi}{4n} \right)^s}{\cos \left(\dfrac{(2p-1)\pi}{4n}\right)} 
			\end{equation}
	then $W(s)$ where $m \in \mathbb{N}$ is denoted by $J(s)$ as follows:
			\begin{equation}
				W(m)=\sum_{k=0}^m \frac{(-1)^k}{(m-k)!} \bigg( \frac{\pi}{2} \bigg)^{m-k} J(k)
			\end{equation}
		\end{lem}
\
		\begin{proof} \small
			\begin{align*}
				W(s)&=\frac{1}{\Gamma(s+1)}	 \lim_{n \to \infty} \frac{1}{n} \sum_{p=1}^n  
					\frac{\left( \dfrac{(2p-1)\pi}{4n} \right)^s}{\cos \left(\dfrac{(2p-1)\pi}{4n}\right)}  
					=\frac{1}{m!} \lim_{n \to \infty} \frac{1}{n} \sum_{p=1}^n
					\frac{\left(\dfrac{(2n-(2n-(2p-1)))\pi}{4n} \right)^m}
						{\cos\left(\dfrac{(2n-(2n-(2p-1)))\pi}{4n} \right)}  \\
					&=\frac{1}{m!} \lim_{n \to \infty} \frac{1}{n} \sum_{p=1}^n 
					\frac{\left(\dfrac{\pi}{2}-\dfrac{(2n-(2p-1))\pi}{4n} \right)^m}{\cos\left(\dfrac{\pi}{2}
					-\dfrac{(2n-(2p-1))\pi}{4n} \right)}  
					=\frac{1}{m!} \lim_{n \to \infty} \frac{1}{n} \sum_{p=1}^n 
					\frac{\left(\dfrac{\pi}{2}-\dfrac{(2n-(2p-1))\pi}{4n} \right)^m}
						{\sin\left(\dfrac{(2n-(2p-1))\pi}{4n} \right)} 
			\end{align*}
	Since $\sum_{p=1}^n f(2n-(2p-1))=\sum_{p=1}^n f((2p-1))$ where $f$ be a real-valued function,
			\begin{align*} 
				W(m)&=\frac{1}{m!} \lim_{n \to \infty} \frac{1}{n} \sum_{p=1}^n  
					\frac{\left(\dfrac{\pi}{2}-\dfrac{(2p-1)\pi}{4n} \right)^m}
						{\sin\left(\dfrac{(2p-1)\pi}{4n} \right)} 
					=\frac{1}{m!} \lim_{n \to \infty} \frac{1}{n} \sum_{p=1}^n \left[
						\frac{\displaystyle \sum_{k=0}^m \left\{ \binom{m}{k} \left(\dfrac{\pi}{2}\right)^{m-k} 
						\left(-\dfrac{(2p-1)\pi}{4n} \right)^k \right\}}{\sin \left( \dfrac{(2p-1)\pi}{4n} \right)} \right] \\
				 	&=\sum_{k=0}^m \frac{(-1)^k}{(m-k)!} \left( \frac{\pi}{2} \right)^{m-k} \left[
						\frac{1}{k!} \lim_{n \to \infty} \frac{1}{n} \sum_{p=1}^n 
						\dfrac{\left(\dfrac{(2p-1)\pi}{4n} \right)^k}{\sin \left( \dfrac{(2p-1)\pi}{4n} \right)} \right] 
					=\sum_{k=0}^m \frac{(-1)^k}{(m-k)!} \left( \frac{\pi}{2} \right)^{m-k} J(k)
			\end{align*}
		\end{proof}

\section{Proof of the Theorems}
$\newline$
	The expression $x(\pi-x)$ where $(0\leq x\leq \pi)$ can be expanded to a Fourier sine series as follows:
		\begin{equation}
			x(\pi-x)=\frac{8}{\pi} \left\{\frac{\sin(x)}{1^3}+\frac{\sin(3x)}{3^3}+\frac{\sin(5x)}{5^3} \cdots \right\}
		\end{equation}
	Using the Dirichlet lambda and beta function values, we have
		\begin{equation} \label{eq:a2}
			\sum_{k=1}^{\infty} \frac{\sin((2k-1)x)}{(2k-1)^3}=\lambda(2) x - \beta(1) \frac{x^2}{2!} 
		\end{equation}
$\newline$
	Let $f_{n}(x)=\sum_{k=1}^{\infty} \dfrac{\sin((2k-1)x)}{(2k-1)^n}$ and $g_{n}(x)=\sum_{k=1}^{\infty} \dfrac{\cos((2k-1)x)}{(2k-1)^n}$, then the multiple integrals on both sides of Eq. (\ref{eq:a2}) with respect to $x$ from 0 to $x$ are given by the functional equations.
		\begin{equation} \label{eq:a3}
			f_{2m+1}(x)=\sum_{k=1}^{\infty} \frac{\sin((2k-1)x)}{(2k-1)^{2m+1}}
				=\sum_{k=1}^{m} \left\{ \lambda(2m-2k+2) \frac{(-1)^{k-1} x^{2k-1}}{(2k-1)!} \right\}
				+ (-1)^{m} \beta(1) \frac{x^{2m}}{(2m)!} 
		\end{equation}
		\begin{equation} \label{eq:a4}
			g_{2m}(x)=\sum_{k=1}^{\infty} \frac{\cos((2k-1)x)}{(2k-1)^{2m}}
				=\sum_{k=1}^{m} \left\{ \lambda(2m-2k+2) \frac{(-1)^{k-1} x^{2k-2}}{(2k-2)!} \right\}
				+ (-1)^{m} \beta(1) \frac{x^{2m-1}}{(2m-1)!} 
		\end{equation}
	where $m \in \mathbb{N}$. The constant of integration is determined by boundary conditions at $f_{n}(0)=0$ and $g_{n}(0)=\lambda(n)$.
$\newline\newline$ 
	If $a_{k}=\sin \left((2k-1) \dfrac{(2p-1)\pi}{4n} \right)$ and $b_{k}=\cos \left((2k-1) \dfrac{(2p-1)\pi}{4n} \right)$ where $p=1, 2, 3, \cdots,n$, periodic sequences $a_k$ and $b_k$ satisfy as follow:
		\begin{align}
			a_k&=(-1)^{m+1} a_{2mn-(k-1)}=(-1)^{m} a_{2mn+k} \\
			b_k&=(-1)^{m} b_{2mn-(k-1)}=(-1)^{m} a_{2mn+k} 
		\end{align}
	where $(1 \leq k \leq n)$ and $k, m \in \mathbb{N}$. For example, if $n=10$ and $k=6$, then $a_6=a_{15}=-a_{26}=-a_{35}=a_{46}= \cdots$ and $b_6=-b_{15}=-b_{26}=b_{35}=b_{46}= \cdots$.
$\newline\newline$ 
	Thus, $f_{2m+1} \left( \dfrac{(2p-1)\pi}{4n} \right)$ and $ g_{2m} \left( \dfrac{(2p-1)\pi}{4n} \right)$ are given by the functional equations.
		\begin{align}
			f_{2m+1} \left( \dfrac{(2p-1)\pi}{4n} \right)
					&=\sum_{q=1}^{n} \left[ \sin\left( \dfrac{(2q-1)(2p-1)\pi}{4n} \right) 
					\sum_{k=1}^{\infty} 
					\frac{(-1)^{k-1}}{ \left\{ (2k-1)2n-(2n-(2q-1)) \right\}^{2m+1}} \right] \nonumber \\
					&+\sum_{q=1}^{n} \left[ \sin\left( \dfrac{(2q-1)(2p-1)\pi}{4n} \right) 
					\sum_{k=1}^{\infty} 
					\frac{(-1)^{k-1}}{ \left\{ (2k-1)2n+(2n-(2q-1)) \right\}^{2m+1}} \right]
		\end{align}
		\begin{align}
			g_{2m} \left( \dfrac{(2p-1)\pi}{4n} \right)
					&=\sum_{q=1}^{n} \left[ \cos\left( \dfrac{(2q-1)(2p-1)\pi}{4n} \right) 
					\sum_{k=1}^{\infty} 
					\frac{(-1)^{k-1}}{ \left\{ (2k-1)2n-(2n-(2q-1)) \right\}^{2m}} \right] \nonumber \\
					&-\sum_{q=1}^{n} \left[ \cos\left( \dfrac{(2q-1)(2p-1)\pi}{4n} \right) 
					\sum_{k=1}^{\infty} 
					\frac{(-1)^{k-1}}{ \left\{ (2k-1)2n+(2n-(2q-1)) \right\}^{2m}} \right]
		\end{align}
$\newline$ 
	When $p$ has the values $1, 2, \cdots, n$, we get $n$ functional equations which can be written as
		\begin{equation} \label{eq:a5}
			\mathbf{F=AX}
		\end{equation}
		\begin{equation}  \label{eq:a6}
			\mathbf{G=BY}
		\end{equation}
	where
		\begin{align}
			\mathbf{A} =
			\left( \begin{array}{cccc}
			\sin \left( \dfrac{\pi}{4n} \right)  & \sin \left( \dfrac{3\pi}{4n} \right) & \ldots & \sin \left( \dfrac{(2n-1)\pi}{4n} \right) \\
			\sin \left( \dfrac{3\pi}{4n} \right) & \sin \left( \dfrac{9\pi}{4n} \right) & \ldots& \sin \left( \dfrac{3(2n-1)\pi}{4n} \right) \\
			\vdots & \vdots & \ddots & \vdots \\
			\sin \left( \dfrac{(2n-1)\pi}{4n} \right) & \sin \left( \dfrac{3(2n-1)\pi}{4n} \right) & \ldots& \sin \left( \dfrac{(2n-1)^2\pi}{4n} \right) \\
			\end{array} \right)
		\end{align}

		\begin{align}
			\mathbf{B} =
			\left( \begin{array}{cccc}
			\cos \left( \dfrac{\pi}{4n} \right)  & \cos \left( \dfrac{3\pi}{4n} \right) & \ldots & \cos \left( \dfrac{(2n-1)\pi}{4n} \right) \\
			\cos \left( \dfrac{3\pi}{4n} \right) & \cos \left( \dfrac{9\pi}{4n} \right) & \ldots& \cos \left( \dfrac{3(2n-1)\pi}{4n} \right) \\
			\vdots & \vdots & \ddots & \vdots \\
			\cos \left( \dfrac{(2n-1)\pi}{4n} \right) & \cos \left( \dfrac{3(2n-1)\pi}{4n} \right) & \ldots& \cos \left( \dfrac{(2n-1)^2\pi}{4n} \right) \\
			\end{array} \right)
		\end{align}

		\begin{align}
			\mathbf{X} =
			\left( \begin{array}{c}
			\displaystyle \sum_{k=1}^{\infty} \dfrac{(-1)^{k-1}}{ \left\{ (2k-1)2n-(2n-1) \right\}^{2m+1}}
			+\displaystyle \sum_{k=1}^{\infty} \dfrac{(-1)^{k-1}}{ \left\{ (2k-1)2n+(2n-1) \right\}^{2m+1}}  \\
			\displaystyle \sum_{k=1}^{\infty} \dfrac{(-1)^{k-1}}{ \left\{ (2k-1)2n-(2n-3) \right\}^{2m+1}}
			+\displaystyle \sum_{k=1}^{\infty} \dfrac{(-1)^{k-1}}{ \left\{ (2k-1)2n+(2n-3) \right\}^{2m+1}}  \\
			\vdots \\
			\displaystyle \sum_{k=1}^{\infty} \dfrac{(-1)^{k-1}}{ \left\{ (2k-1)2n-(2n-(2n-1)) \right\}^{2m+1}}
			+\displaystyle \sum_{k=1}^{\infty} \dfrac{(-1)^{k-1}}{ \left\{ (2k-1)2n+(2n-(2n-1)) \right\}^{2m+1}}  \\
			\end{array} \right)
		\end{align}

		\begin{align}
			\mathbf{Y} =
			\left( \begin{array}{c}
			\displaystyle \sum_{k=1}^{\infty} \dfrac{(-1)^{k-1}}{ \left\{ (2k-1)2n-(2n-1) \right\}^{2m}}
			-\displaystyle \sum_{k=1}^{\infty} \dfrac{(-1)^{k-1}}{ \left\{ (2k-1)2n+(2n-1) \right\}^{2m}}  \\
			\displaystyle \sum_{k=1}^{\infty} \dfrac{(-1)^{k-1}}{ \left\{ (2k-1)2n-(2n-3) \right\}^{2m}}
			-\displaystyle \sum_{k=1}^{\infty} \dfrac{(-1)^{k-1}}{ \left\{ (2k-1)2n+(2n-3) \right\}^{2m}}  \\
			\vdots \\
			\displaystyle \sum_{k=1}^{\infty} \dfrac{(-1)^{k-1}}{ \left\{ (2k-1)2n-(2n-(2n-1)) \right\}^{2m}}
			-\displaystyle \sum_{k=1}^{\infty} \dfrac{(-1)^{k-1}}{ \left\{ (2k-1)2n+(2n-(2n-1)) \right\}^{2m}}  \\
			\end{array} \right)
		\end{align}

		\begin{align}
			\mathbf{F} =
			\left( \begin{array}{c}
			f_{2m+1} \left( \dfrac{\pi}{4n} \right) \\
			f_{2m+1} \left( \dfrac{3\pi}{4n} \right) \\
			\vdots \\
			f_{2m+1} \left( \dfrac{(2n-1)\pi}{4n} \right) \\
			\end{array} \right)
			 \qquad 
			\mathbf{G} =
			\left( \begin{array}{c}
			g_{2m} \left( \dfrac{\pi}{4n} \right) \\
			g_{2m} \left( \dfrac{3\pi}{4n} \right) \\
			\vdots \\
			g_{2m} \left( \dfrac{(2n-1)\pi}{4n} \right) \\
			\end{array} \right)
		\end{align}
$\newline$ 
	To calculate $\mathbf{X}$ and $\mathbf{Y}$ (Eq. (\ref{eq:a5}) and Eq. (\ref{eq:a6})), we need to the Lemma $\ref{lem:3}$ and $\ref{lem:4}$. By the Lemma $\ref{lem:3}$ and $\ref{lem:4}$, the $\mathbf{X=A^{-1}F}=\dfrac{2}{n}\mathbf{AF}$ and $\mathbf{Y=B^{-1}G=}\dfrac{2}{n}\mathbf{BG}$ as follows:
		\begin{align*} 
			\mathbf{X} =
			\left( \begin{array}{c}
			\displaystyle \sum_{k=1}^{\infty} \dfrac{(-1)^{k-1}}{ \left\{ (2k-1)2n-(2n-1) \right\}^{2m+1}}
			+\displaystyle \sum_{k=1}^{\infty} \dfrac{(-1)^{k-1}}{ \left\{ (2k-1)2n+(2n-1) \right\}^{2m+1}}  \\
			\displaystyle \sum_{k=1}^{\infty} \dfrac{(-1)^{k-1}}{ \left\{ (2k-1)2n-(2n-3) \right\}^{2m+1}}
			+\displaystyle \sum_{k=1}^{\infty} \dfrac{(-1)^{k-1}}{ \left\{ (2k-1)2n+(2n-3) \right\}^{2m+1}}  \\
			\vdots \\
			\displaystyle \sum_{k=1}^{\infty} \dfrac{(-1)^{k-1}}{ \left\{ (2k-1)2n-(2n-(2n-1)) \right\}^{2m+1}}
			+\displaystyle \sum_{k=1}^{\infty} \dfrac{(-1)^{k-1}}{ \left\{ (2k-1)2n+(2n-(2n-1)) \right\}^{2m+1}}  \\
			\end{array} \right) 
		\end{align*}
		\begin{align}  \label{eq:a7} \small
			=\dfrac{2}{n}\left( \begin{array}{cccc}
			\sin \left( \dfrac{\pi}{4n} \right)  & \sin \left( \dfrac{3\pi}{4n} \right) & \ldots & \sin \left( \dfrac{(2n-1)\pi}{4n} \right) \\
			\sin \left( \dfrac{3\pi}{4n} \right) & \sin \left( \dfrac{9\pi}{4n} \right) & \ldots& \sin \left( \dfrac{3(2n-1)\pi}{4n} \right) \\
			\vdots & \vdots & \ddots & \vdots \\
			\sin \left( \dfrac{(2n-1)\pi}{4n} \right) & \sin \left( \dfrac{3(2n-1)\pi}{4n} \right) & \ldots& \sin \left( \dfrac{(2n-1)^2\pi}{4n} \right) \\
			\end{array} \right)
			\left( \begin{array}{c}
			f_{2m+1} \left( \dfrac{\pi}{4n} \right) \\
			f_{2m+1} \left( \dfrac{3\pi}{4n} \right) \\
			\vdots \\
			f_{2m+1} \left( \dfrac{(2n-1)\pi}{4n} \right) \\
			\end{array} \right)
		\end{align}
$\newline$ 
		\begin{align*}
			\mathbf{Y} =
			\left( \begin{array}{c}
			\displaystyle \sum_{k=1}^{\infty} \dfrac{(-1)^{k-1}}{ \left\{ (2k-1)2n-(2n-1) \right\}^{2m}}
			-\displaystyle \sum_{k=1}^{\infty} \dfrac{(-1)^{k-1}}{ \left\{ (2k-1)2n+(2n-1) \right\}^{2m}}  \\
			\displaystyle \sum_{k=1}^{\infty} \dfrac{(-1)^{k-1}}{ \left\{ (2k-1)2n-(2n-3) \right\}^{2m}}
			-\displaystyle \sum_{k=1}^{\infty} \dfrac{(-1)^{k-1}}{ \left\{ (2k-1)2n+(2n-3) \right\}^{2m}}  \\
			\vdots \\
			\displaystyle \sum_{k=1}^{\infty} \dfrac{(-1)^{k-1}}{ \left\{ (2k-1)2n-(2n-(2n-1)) \right\}^{2m}}
			-\displaystyle \sum_{k=1}^{\infty} \dfrac{(-1)^{k-1}}{ \left\{ (2k-1)2n+(2n-(2n-1)) \right\}^{2m}}  \\
			\end{array} \right)
		\end{align*}
		\begin{align}  \label{eq:a8} \small
			=\dfrac{2}{n} \left( \begin{array}{cccc}
			\cos \left( \dfrac{\pi}{4n} \right)  & \cos \left( \dfrac{3\pi}{4n} \right) & \ldots & \cos \left( \dfrac{(2n-1)\pi}{4n} \right) \\
			\cos \left( \dfrac{3\pi}{4n} \right) & \cos \left( \dfrac{9\pi}{4n} \right) & \ldots& \cos \left( \dfrac{3(2n-1)\pi}{4n} \right) \\
			\vdots & \vdots & \ddots & \vdots \\
			\cos \left( \dfrac{(2n-1)\pi}{4n} \right) & \cos \left( \dfrac{3(2n-1)\pi}{4n} \right) & \ldots& \cos \left( \dfrac{(2n-1)^2\pi}{4n} \right) \\
			\end{array} \right) 
			\left( \begin{array}{c}
			g_{2m} \left( \dfrac{\pi}{4n} \right) \\
			g_{2m} \left( \dfrac{3\pi}{4n} \right) \\
			\vdots \\
			g_{2m} \left( \dfrac{(2n-1)\pi}{4n} \right) \\
			\end{array} \right)
		\end{align}
\\
\subsection{Dirichlet Lambda Function at Odd Positive Integers} \
$\newline\newline$
	In Eq. (\ref{eq:a7}), we know that sum of all elements in a matrix $\mathbf{X}$ is equal to $\lambda(2m+1)$
		\begin{align*}
			\lambda(2m+1)=\sum_{k=1}^{\infty} \frac{1}{(2k-1)^{2m+1}}=\lim_{n \to \infty} \sum_{k=1}^{n} \mathbf{X}_{k, 1}
		\end{align*}
	where
		\begin{align*} \small
			\mathbf{X} 
			=\dfrac{2}{n}\left( \begin{array}{cccc}
			\sin \left( \dfrac{\pi}{4n} \right)  & \sin \left( \dfrac{3\pi}{4n} \right) & \ldots & \sin \left( \dfrac{(2n-1)\pi}{4n} \right) \\
			\sin \left( \dfrac{3\pi}{4n} \right) & \sin \left( \dfrac{9\pi}{4n} \right) & \ldots& \sin \left( \dfrac{3(2n-1)\pi}{4n} \right) \\
			\vdots & \vdots & \ddots & \vdots \\
			\sin \left( \dfrac{(2n-1)\pi}{4n} \right) & \sin \left( \dfrac{3(2n-1)\pi}{4n} \right) & \ldots& \sin \left( \dfrac{(2n-1)^2\pi}{4n} \right) \\
			\end{array} \right)
			\left( \begin{array}{c}
			f_{2m+1} \left( \dfrac{\pi}{4n} \right) \\
			f_{2m+1} \left( \dfrac{3\pi}{4n} \right) \\
			\vdots \\
			f_{2m+1} \left( \dfrac{(2n-1)\pi}{4n} \right) \\
			\end{array} \right)
		\end{align*}
		\begin{align*} \small
			=\dfrac{2}{n} \left( \begin{array}{cccc}
			\left( \begin{array}{c}
			\sin \left( \dfrac{\pi}{4n} \right) \\
			\sin \left( \dfrac{3\pi}{4n} \right) \\
			\vdots \\
			\sin \left( \dfrac{(2n-1)\pi}{4n} \right) \\
			\end{array} \right) f_{2m+1} \left( \dfrac{\pi}{4n}\right) + \cdots +
			\left( \begin{array}{c}
			\sin \left( \dfrac{(2n-1)\pi}{4n} \right) \\
			\sin \left( \dfrac{3(2n-1)\pi}{4n} \right) \\
			\vdots \\
			\sin \left( \dfrac{(2n-1)^2\pi}{4n} \right) \\
			\end{array} \right) f_{2m+1} \left( \dfrac{(2n-1)\pi}{4n}\right)
			\end{array} \right)
		\end{align*}
	Thus, sum of all elements in a matrix $\mathbf{X}$ is represented as
		\begin{align*}
			\lambda(2m+1)= \lim_{n \to \infty} \frac{2}{n} \sum_{p=1}^{n}
				\left[ f_{2m+1} \left( \frac{(2p-1)\pi}{4n} \right) 
				\sum_{q=1}^{n} \sin \left( \frac{(2p-1)(2q-1)\pi}{4n} \right) \right]
		\end{align*}
	By using the Lemma $\ref{lem:1}$, we have
		\begin{align*}
			\lambda(2m+1)&= \lim_{n \to \infty} \frac{2}{n} \sum_{p=1}^{n}
				\left[ f_{2m+1} \left( \frac{(2p-1)\pi}{4n} \right) 
				\frac{\sin^2 \left( \dfrac{(2p-1)\pi}{4} \right)}{\sin \left( \dfrac{(2p-1)\pi}{4n} \right)} \right] \\
				&= \lim_{n \to \infty} \frac{1}{n} \sum_{p=1}^{n}
				\left[ f_{2m+1} \left( \frac{(2p-1)\pi}{4n} \right) 
				\frac{1}{\sin \left( \dfrac{(2p-1)\pi}{4n} \right)} \right]
		\end{align*}
	Using the Eq. (\ref{eq:a3}), we have
		\begin{align*}
			\lambda(2m+1)&= \lim_{n \to \infty} \frac{1}{n} \sum_{p=1}^{n} \left[
			\sum_{k=1}^{m} \dfrac{(-1)^{k-1} \lambda(2m-2k+2)}{(2k-1)!} 
			\dfrac{ \left( \dfrac{(2p-1)\pi}{4n} \right)^{2k-1}}{\sin \left( \dfrac{(2p-1)\pi}{4n} \right)}  \right] \\
			&+ \lim_{n \to \infty} \frac{1}{n} \sum_{p=1}^{n} \left[
			\frac{(-1)^{m} \beta(1)}{(2m)!} \frac{\left( \dfrac{(2p-1)\pi}{4n} \right)^{2m}}{\sin \left( \dfrac{(2p-1)\pi}{4n} \right)} \right]
		\end{align*}
	Using the Lemma $\ref{lem:5}$, we have
		\begin{align*}
			\lambda(2m+1)&=\sum_{k=1}^{m} \left[ (-1)^{k-1} \lambda(2m-2k+2) J(2k-1) \right] 
			+(-1)^{m} \beta(1) J(2m) \\
		\end{align*}
	The proof of Theorem 1 was completed.
\\

\subsection{Dirichlet Beta Function at Even Positive Integers} \
$\newline\newline$
	In Eq. (\ref{eq:a8}), we know that sum of all elements in a matrix $\mathbf{Y}$ is equal to $\lambda(2m)$
		\begin{align*}
			\lambda(2m)=\sum_{k=1}^{\infty} \frac{1}{(2k-1)^{2m}}=\lim_{n \to \infty} \sum_{k=1}^{n} \mathbf{Y}_{k, 1}
		\end{align*}
	where
		\begin{align*}
			\mathbf{Y} =
			\dfrac{2}{n} \left( \begin{array}{cccc}
			\cos \left( \dfrac{\pi}{4n} \right)  & \cos \left( \dfrac{3\pi}{4n} \right) & \ldots & \cos \left( \dfrac{(2n-1)\pi}{4n} \right) \\
			\cos \left( \dfrac{3\pi}{4n} \right) & \cos \left( \dfrac{9\pi}{4n} \right) & \ldots& \cos \left( \dfrac{3(2n-1)\pi}{4n} \right) \\
			\vdots & \vdots & \ddots & \vdots \\
			\cos \left( \dfrac{(2n-1)\pi}{4n} \right) & \cos \left( \dfrac{3(2n-1)\pi}{4n} \right) & \ldots& \cos \left( \dfrac{(2n-1)^2\pi}{4n} \right) \\
			\end{array} \right) 
			\left( \begin{array}{c}
			g_{2m} \left( \dfrac{\pi}{4n} \right) \\
			g_{2m} \left( \dfrac{3\pi}{4n} \right) \\
			\vdots \\
			g_{2m} \left( \dfrac{(2n-1)\pi}{4n} \right) \\
			\end{array} \right)
		\end{align*}
	In order to obtain the expression $\beta(2m)$, we define the matrix $\mathbf{Z}$ as follows:
		\begin{align*} \small
			\mathbf{Z} =
			\dfrac{2}{n} \left( \begin{array}{cccc}
			\left( \begin{array}{c}
			\cos \left( \dfrac{\pi}{4n} \right) \\
			-\cos \left( \dfrac{3\pi}{4n} \right) \\
			\vdots \\
			(-1)^{n-1}\cos \left( \dfrac{(2n-1)\pi}{4n} \right) \\
			\end{array} \right) g_{2m} \left( \dfrac{\pi}{4n}\right) + \cdots +
			\left( \begin{array}{c}
			\cos \left( \dfrac{(2n-1)\pi}{4n} \right) \\
			-\cos \left( \dfrac{3(2n-1)\pi}{4n} \right) \\
			\vdots \\
			(-1)^{n-1}\cos \left( \dfrac{(2n-1)^2\pi}{4n} \right) \\
			\end{array} \right) g_{2m} \left( \dfrac{(2n-1)\pi}{4n}\right)
			\end{array} \right)
		\end{align*}
	Then sum of all elements in a matrix $\mathbf{Z}$ is equal to $\beta(2m)$.
		\begin{align*} 
			\beta(2m)=\lim_{n \to \infty} \frac{2}{n} \sum_{p=1}^{n} \left[
			g_{2m} \left( \frac{(2p-1)\pi}{4n} \right) \sum_{q=1}^{n} (-1)^{q-1} \cos \left( \frac{(2p-1)(2q-1)\pi}{4n} \right) \right]
		\end{align*}
	By using the Lemma $\ref{lem:2}$, we have
		\begin{align*}
			\beta(2m)&= \lim_{n \to \infty} \frac{2}{n} \sum_{p=1}^{n}
				\left[ g_{2m} \left( \frac{(2p-1)\pi}{4n} \right) 
				\frac{\sin^2 \left( \dfrac{n}{2} \left( \pi- \dfrac{(2p-1)\pi}{2n} \right) \right)}{2\cos \left( \dfrac{(2p-1)\pi}{4n} \right)} \right] \\
				&= \lim_{n \to \infty} \frac{1}{n} \sum_{p=1}^{n}
				\left[ g_{2m} \left( \frac{(2p-1)\pi}{4n} \right) 
				\frac{1}{\cos \left( \dfrac{(2p-1)\pi}{4n} \right)} \right] \\
		\end{align*}
	Using the Eq. (\ref{eq:a4}), we have
		\begin{align*}
			\beta(2m)&= \lim_{n \to \infty} \frac{1}{n} \sum_{p=1}^{n} \left[
			\sum_{k=1}^{m} \dfrac{(-1)^{k-1} \lambda(2m-2k+2)}{(2k-2)!} 
			\dfrac{ \left( \dfrac{(2p-1)\pi}{4n} \right)^{2k-2}}{\cos \left( \dfrac{(2p-1)\pi}{4n} \right)}  \right] \\
			&+ \lim_{n \to \infty} \frac{1}{n} \sum_{p=1}^{n} \left[
			\frac{(-1)^{m} \beta(1)}{(2m-1)!} \frac{\left( \dfrac{(2p-1)\pi}{4n} \right)^{2m-1}}{\cos \left( \dfrac{(2p-1)\pi}{4n} \right)} \right]
		\end{align*}
	Using the Lemma $\ref{lem:6}$, we have
		\begin{align}
			\beta(2m)&=\sum_{k=1}^{m} \left[ (-1)^{k-1} \lambda(2m-2k+2) W(2k-2) \right] + (-1)^{m} \beta(1) W(2m-1) \\
					&=\sum_{k=1}^{m} \left[ (-1)^{k-1} \lambda(2m-2k+2) \sum_{q=0}^{2k-2} 
						\left\{ (-1)^{q} \frac{1}{\left\{(2k-2)-q \right\}!} 
						\left( \frac{\pi}{2} \right)^{(2k-2)-q} J(q) \right\} \right] \nonumber \\
					&+(-1)^{m} \beta(1) \sum_{q=0}^{2m-1} \left\{ (-1)^{q} \frac{1}{\left\{(2m-1)-q \right\}!}
						\left( \frac{\pi}{2} \right)^{(2m-1)-q} J(q) \right\} \nonumber
		\end{align}
	The index of summation $q$ takes on integer values from 0 to $2k-2$. 
\newline\newline
	Now, expand the inner summation(which involves $q$).

		\begin{align*}
			\beta(2m)&=\sum_{k=1}^{m} \left[ (-1)^{k-1} \frac{\lambda(2m-2k+2)}{(2k-2)!} \left( \frac{\pi}{2} \right )^{2k-2} J(+0) \right]
					+\frac{(-1)^{m} \beta(1)}{(2m-1)!} \left( \frac{\pi}{2} \right)^{2m-1} J(+0) \\
					&+\sum_{k=2}^{m} \left[ (-1)^{k} \frac{\lambda(2m-2k+2)}{\{(2k-2)-1\}!} \left( \frac{\pi}{2} \right )^{(2k-2)-1} J(1) \right]
					+\frac{(-1)^{m+1} \beta(1)}{\{(2m-1)-1\}!} \left( \frac{\pi}{2} \right)^{(2m-1)-1} J(1) \\
					&+\sum_{k=2}^{m} \left[ (-1)^{k+1} \frac{\lambda(2m-2k+2)}{\{(2k-2)-2\}!} \left( \frac{\pi}{2} \right )^{(2k-2)-2} J(2) \right]
					+\frac{(-1)^{m+2} \beta(1)}{\{(2m-1)-2\}!} \left( \frac{\pi}{2} \right)^{(2m-1)-2} J(2) + \cdots \\
					&+\sum_{k=m}^{m} \left[ (-1)^{k+2m-3} \frac{\lambda(2m-2k+2)}{\{(2k-2)-(2m-2)\}!} \left( \frac{\pi}{2} \right )^{(2k-2)-(2m-2)} J(2k-2) \right] \\
					&+\frac{(-1)^{3m-1} \beta(1)}{\{(2m-1)-(2m-1)\}!} \left( \frac{\pi}{2} \right)^{(2m-1)-(2m-1)} J(2m-1)
		\end{align*}
\newline

	Change the index of summation $k$ so that it would start from 1.
		\begin{align*}
			\beta(2m)&=\left[ \sum_{k=1}^{m} \left\{ (-1)^{k-1} \frac{\lambda(2m-2k+2)}{(2k-2)!} \left( \frac{\pi}{2} \right )^{2k-2} \right\}
					+\frac{(-1)^{m} \beta(1)}{(2m-1)!} \left( \frac{\pi}{2} \right)^{2m-1} \right] J(+0) \\
					&+\left[ \sum_{k=1}^{m-1} \left\{ (-1)^{k-1} \frac{\lambda(2m-2k)}{(2k-1)!} \left( \frac{\pi}{2} \right )^{2k-1} \right\}
					+\frac{(-1)^{m-1} \beta(1)}{(2m-2)!} \left( \frac{\pi}{2} \right)^{2m-2} \right] J(1) \\
					&+\left[ \sum_{k=1}^{m-1} \left\{ (-1)^{k} \frac{\lambda(2m-2k)}{(2k-2)!} \left( \frac{\pi}{2} \right )^{2k-2} \right\}
					+\frac{(-1)^{m} \beta(1)}{(2m-3)!} \left( \frac{\pi}{2} \right)^{2m-3} \right] J(2) + \cdots \\
					&+\left[ (-1)^{m} \lambda(2) J(2m-2) \right] +(-1)^{m-1} \beta(1) J(2m-1) 
		\end{align*}
	Using the Eq. (\ref{eq:a3}) and Eq. (\ref{eq:a4})
		\begin{align*}
			\beta(2m)&=\sum_{k=1}^{m} \left[ (-1)^{k-1} g_{2m-2k+2} \left( \frac{\pi}{2} \right) J(2k-2)+(-1)^{k-1} f_{2m-2k+1} \left( \frac{\pi}{2} \right) J(2k-1) \right]
		\end{align*}
	Since $g_{2m} \left( \dfrac{\pi}{2} \right)=0$ and $f_{2m+1} \left( \dfrac{\pi}{2} \right)=\beta(2m-1)$ (See  Eq. $(\ref{eq:a3})$ and Eq. $(\ref{eq:a4})$)
		\begin{align*}
			\beta(2m)=\sum_{k=1}^{m} \left[ (-1)^{k-1} \beta(2m-2k+1) J(2k-1) \right]
		\end{align*}
	The proof of Theorem 2 was completed.
\\

\section{THE INTEGRAL FUNCTION $J(n)$}

%LEMMA 7

		\begin{lem} \label{lem:7}
	The function $-\dfrac{1}{2} \ln \left( \tan \dfrac{x}{2} \right)$ can be expanded as an infinite series,
		\begin{align}
			\sum_{k=1}^{\infty} \frac{\cos((2k-1)x)}{2k-1}=-\dfrac{1}{2} \ln \left( \tan \dfrac{x}{2} \right) 
		\end{align}
	where $x \in \mathbb{R}$
		\end{lem}
\
		\begin{proof}
	Let $f(x)=\sum_{k=1}^{\infty} e^{i(2k-1)x}$, then we have
		\begin{align*}
			f(x)=\sum_{k=1}^{\infty} e^{i(2k-1)x}=\frac{e^{ix}}{1-e^{2ix}}=\frac{1}{e^{-ix}-e^{ix}}=\frac{1}{2i} \frac{2i}{e^{-ix}-e^{ix}}=\frac{i}{2} \csc(x)
		\end{align*}
	By integrating the $f(x)$, we have
		\begin{align*}
			&\frac{1}{i} e^{ix} + \frac{1}{3i} e^{3ix} + \frac{1}{5i} e^{5ix} + \cdots = \frac{i}{2} \ln \left( \tan \frac{x}{2} \right)+C \\
			&e^{ix} + \frac{1}{3} e^{3ix} + \frac{1}{5} e^{5ix} + \cdots = -\frac{1}{2} \ln \left( \tan \frac{x}{2} \right)+Ci 
		\end{align*}
	where $C$ is constant of integration.
\newline\newline
	Taking the real part,
		\begin{align*}
			\cos(x) + \frac{1}{3} \cos(3x) + \frac{1}{5} \cos(5x) + \cdots = -\frac{1}{2} \ln \left( \tan \frac{x}{2} \right) 
		\end{align*}
		\end{proof}

%LEMMA 8

		\begin{lem} \label{lem:8}
	The Euler number $E_{n}$ is represented as
		\begin{align}
			\frac{d^{2n}}{dx^{2n}} \csc\left( \frac{\pi}{2} \right)=(-1)^{n} E_{2n}
		\end{align}
	where $n \in \left\{ \mathbb{N}, 0 \right\}$
		\end{lem}
\
		\begin{proof}
	The expression for $\csc(x)$ can be expanded to a Taylor series at $x=\pi/2$ as follows:
		\begin{align*}
			\csc(x)=1+\frac{1}{2!}\left(x-\frac{\pi}{2}\right)^{2}+\frac{5}{4!}\left(x-\frac{\pi}{2}\right)^{4}+\frac{61}{6!}\left(x-\frac{\pi}{2}\right)^{6}
				+\cdots = \sum_{m=0}^{\infty} \frac{\left|E_{m}\right|}{m!} \left(x-\frac{\pi}{2} \right)^{m}
		\end{align*}
	The definition of the $m-th$ term of a Taylor series at $x=\pi/2$ is
		\begin{align*}
			\left\{ \frac{d^{m}}{dx^{m}} f\left( \frac{\pi}{2} \right) \right\} \frac{1}{m!} \left(x-\frac{\pi}{2} \right)^{m} 
		\end{align*}
	If $m=2n$, then $E_{m}=E_{2n}$ and If $m=2n+1$, then $E_{2m}=0$. 
	\newline\newline
	Therefore,
		\begin{align*}
			\frac{d^{2n}}{dx^{2n}} \csc\left( \frac{\pi}{2} \right)=(-1)^{n} E_{2n}
		\end{align*}
		\end{proof}

%THEOREM 3

		\begin{thm} \label{thm:3}
	The function $J(n)$ where $n \in \mathbb{N}$ can be expressed as an infinite series,
		\begin{equation}
			J(n)=\sum_{k=0}^{\infty} \frac{(-1)^{k} E_{2k}}{(n+2k+1)!} \left( \frac{\pi}{2} \right)^{n+2k}
		\end{equation} 
	where $E_{k}$ is Euler number.
		\end{thm}
\
		\begin{proof}
	The function $J(n)$ where $n \in \mathbb{N}$ is defined as
		\begin{align*}
			J(n)=\frac{1}{n!} \frac{2}{\pi} \int_{0}^{\frac{\pi}{2}} \frac{x^n}{\sin(x)} dx
		\end{align*}
	Integrating by parts,
		\begin{align*}
			J(n)&=\frac{1}{n!} \frac{2}{\pi} \left[ \frac{x^{n+1}\left\{\csc(x)\right\}}{(n+1)} - \frac{x^{n+2}\left\{ \dfrac{d}{dx} \csc(x)\right\}}{(n+1)(n+2)}
				+ \frac{x^{n+3}\left\{ \dfrac{d^2}{dx^2} \csc(x)\right\}}{(n+1)(n+2)(n+3)} - \cdots \right]^{\frac{\pi}{2}} _{0} \\
				&=\frac{2}{\pi} \left[ \frac{x^{n+1}}{(n+1)!}\left\{\csc(x)\right\} - \frac{x^{n+2}}{(n+2)!} \left\{ \dfrac{d}{dx} \csc(x)\right\}
				+ \frac{x^{n+3}}{(n+3)!}\left\{ \dfrac{d^2}{dx^2} \csc(x)\right\} - \cdots \right]^{\frac{\pi}{2}} _{0}
		\end{align*}
	By using Lemma $\ref{lem:8}$, we have
		\begin{align*}
			J(n)&= \frac{2}{\pi} \left[ \frac{E_0}{(n+1)!} \left( \frac{\pi}{2} \right)^{n+1} + \frac{-E_2}{(n+3)!} \left( \frac{\pi}{2} \right)^{n+3}
				+ \frac{E_4}{(n+5)!} \left( \frac{\pi}{2} \right)^{n+5} +\cdots \right] \\
				&=\left[ \frac{E_0}{(n+1)!} \left( \frac{\pi}{2} \right)^{n} - \frac{E_2}{(n+3)!} \left( \frac{\pi}{2} \right)^{n+2}
				+ \frac{E_4}{(n+5)!} \left( \frac{\pi}{2} \right)^{n+4} -\cdots \right] \\
				&=\sum_{k=0}^{\infty} \frac{(-1)^{k} E_{2k}}{(n+2k+1)!} \left( \frac{\pi}{2} \right)^{n+2k}
		\end{align*}
		The proof of Theorem 3 was completed.
		\end{proof}

$\newline$
%THEOREM 4

		\begin{thm} \label{thm:4}
	The function $J(2n-1)$ and $J(2n)$ where $n \in \mathbb{N}$ can be calculated directly in special forms as
		\begin{gather}
			\frac{\pi}{4}J(2n-1)=(-1)^{n-1} \sum_{k=0}^{n-1} \left[ (-1)^{k} \beta(2n-2k) \frac{1}{(2k)!} \left( \frac{\pi}{2} \right)^{2k} \right] \\
			\frac{\pi}{4}J(2n)=(-1)^{n} \left[ \lambda(2n+1)- \sum_{k=0}^{n-1} \left\{ (-1)^{k} \beta(2n-2k) \frac{1}{(2k+1)!} \left( \frac{\pi}{2} \right)^{2k+1} \right\} \right]
		\end{gather} 
	where $E_{k}$ is Euler number.
		\end{thm}
\
		\begin{proof}
	The expression for $\csc(x)$ can be expanded to a Taylor series at $x=\pi/2$ as follows:
		\begin{align*}
			\csc(x)=1+\frac{1}{2!}\left(x-\frac{\pi}{2}\right)^{2}+\frac{5}{4!}\left(x-\frac{\pi}{2}\right)^{4}+\frac{61}{6!}\left(x-\frac{\pi}{2}\right)^{6}
				+\cdots = \sum_{k=0}^{\infty} \frac{(-1)^{k} E_{2k}}{(2k)!} \left(x-\frac{\pi}{2} \right)^{2k}
		\end{align*}
	Integrating both sides of the formula with respect to $x$, we have
		\begin{align*}
			\ln\left(\tan\frac{x}{2}\right) =\left(x-\frac{\pi}{2} \right) + \frac{1}{3!}\left(x-\frac{\pi}{2}\right)^{3}+\frac{5}{5!}\left(x-\frac{\pi}{2}\right)^{5}
				+\cdots = \sum_{k=0}^{\infty} \frac{(-1)^{k} E_{2k}}{(2k+1)!} \left(x-\frac{\pi}{2} \right)^{2k+1}
		\end{align*}
	The constant of integration is determined by $x=\pi/2$. 
\newline
	By using Lemma $\ref{lem:7}$, we have
		\begin{align*}
			\sum_{k=0}^{\infty} \frac{(-1)^{k} E_{2k}}{(2k+1)!} \left(x-\frac{\pi}{2} \right)^{2k+1}
			=2 \sum_{k=0}^{\infty} \left\{ \frac{-\cos((2k+1)x)}{(2k+1)} \right\} 
		\end{align*}
\newline
	The multiple integral on both sides with respect to $x$ is given by the functional equations. The constant of integration is determined by $x=\pi/2$.
		\begin{align*}
			\sum_{k=0}^{\infty}& \frac{(-1)^{k} E_{2k}}{(2k+2n)!} \left(x-\frac{\pi}{2} \right)^{2k+2n} \\
			&=2 (-1)^{n} \sum_{k=0}^{\infty} \left\{ \frac{\sin((2k+1)x)}{(2k+1)^{2n}} \right\} 
				+ 2(-1)^{n-1} \sum_{k=0}^{n-1} \left[ (-1)^{k} \frac{\beta(2n-2k)}{(2k)!} \left( x-\frac{\pi}{2} \right)^{2k} \right] \\
			\sum_{k=0}^{\infty}& \frac{(-1)^{k} E_{2k}}{(2k+2n+1)!} \left(x-\frac{\pi}{2} \right)^{2k+2n+1} \\
			&=2 (-1)^{n} \sum_{k=0}^{\infty} \left\{ \frac{\cos((2k+1)x)}{(2k+1)^{2n+1}} \right\} 
				+ 2(-1)^{n-1} \sum_{k=0}^{n-1} \left[ (-1)^{k} \frac{\beta(2n-2k)}{(2k+1)!} \left( x-\frac{\pi}{2} \right)^{2k+1} \right] 
		\end{align*}
\newline
	By substituting $x=\pi$, we have
		\begin{align*}
			&\sum_{k=0}^{\infty} \frac{(-1)^{k} E_{2k}}{(2k+2n)!} \left(\frac{\pi}{2} \right)^{2k+2n} 
				=2(-1)^{n-1} \sum_{k=0}^{n-1} \left[ (-1)^{k} \frac{\beta(2n-2k)}{(2k)!} \left(\frac{\pi}{2} \right)^{2k} \right] \\
			&\sum_{k=0}^{\infty} \frac{(-1)^{k} E_{2k}}{(2k+2n+1)!} \left(\frac{\pi}{2} \right)^{2k+2n+1} 
			=2 (-1)^{n} \left[ \lambda(2n+1) - \sum_{k=0}^{n-1} \left\{ (-1)^{k} \frac{\beta(2n-2k)}{(2k+1)!} \left( x-\frac{\pi}{2} \right)^{2k+1} \right\} \right]
		\end{align*}
\newline
	Using the Theorem $\ref{thm:3}$, we have
		\begin{align*}
			&\frac{\pi}{4}J(2n-1)=(-1)^{n-1} \sum_{k=0}^{n-1} \left[ (-1)^{k} \beta(2n-2k) \frac{1}{(2k)!} \left( \frac{\pi}{2} \right)^{2k} \right] \\
			&\frac{\pi}{4}J(2n)=(-1)^{n} \left[ \lambda(2n+1)- \sum_{k=0}^{n-1} \left\{ (-1)^{k} \beta(2n-2k) \frac{1}{(2k+1)!} \left( \frac{\pi}{2} \right)^{2k+1} \right\} \right] 
		\end{align*}
\newline
	The proof of Theorem 4 was completed. \\
		\end{proof}

%REMARK

		\begin{rem} \label{rem:1}
	Similarly to the Theorem $\ref{thm:4}$, the expressions $\dfrac{1}{(2m-1)!} \left( \dfrac{\pi}{2} \right)^{2m-1}$ and $\dfrac{1}{(2m)!} \left( \dfrac{\pi}{2} \right)^{2m}$ where $m \in \mathbb{N}$ can be calculated directly special forms as
		\begin{align}
			&\frac{\pi}{4} \left\{ \dfrac{1}{(2m-1)!} \left( \dfrac{\pi}{2} \right)^{2m-1} \right\}
				=(-1)^{m-1} \sum_{k=0}^{m-1} \left[ (-1)^{k} \lambda(2m-2k) \frac{1}{(2k)!} \left( \frac{\pi}{2} \right)^{2k} \right] \\
			&\frac{\pi}{4}  \left\{ \dfrac{1}{(2m)!} \left( \dfrac{\pi}{2} \right)^{2m} \right\}
				=(-1)^{m} \left[ \beta(2m+1)- \sum_{k=0}^{m-1} \left\{ (-1)^{k} \lambda(2m-2k) \frac{1}{(2k+1)!} \left( \frac{\pi}{2} \right)^{2k+1} \right\} \right]
		\end{align} 
		\end{rem}
\
		\begin{proof}
	By substituting $x=\pi/2$ in Eq. (\ref{eq:a3}) and Eq. (\ref{eq:a4}), we have
		\begin{align*}
			\beta(2m+1)&=\sum_{k=1}^{m} \left\{ \lambda(2m-2k+2) \frac{(-1)^{k-1} {(\pi/2)}^{2k-1}}{(2k-1)!} \right\}
				+ (-1)^{m} \beta(1) \frac{{(\pi/2)}^{2m}}{(2m)!} \\
			0&=\sum_{k=1}^{m} \left\{ \lambda(2m-2k+2) \frac{(-1)^{k-1} {(\pi/2)}^{2k-2}}{(2k-2)!} \right\}
				+ (-1)^{m} \beta(1) \frac{{(\pi/2)}^{2m-1}}{(2m-1)!}  
		\end{align*}
	Change the index of summation $k$ so that it would start from 0.
		\begin{align*}
			\beta(2m+1)&=\sum_{k=0}^{m-1} \left\{ \lambda(2m-2k) \frac{(-1)^{k} {(\pi/2)}^{2k+1}}{(2k+1)!} \right\}
				+ (-1)^{m} \beta(1) \frac{{(\pi/2)}^{2m}}{(2m)!} \\
			0&=\sum_{k=0}^{m-1} \left\{ \lambda(2m-2k) \frac{(-1)^{k-1} {(\pi/2)}^{2k}}{(2k)!} \right\}
				+ (-1)^{m} \beta(1) \frac{{(\pi/2)}^{2m-1}}{(2m-1)!}  
		\end{align*}
	Therefore,
		\begin{align*}
			&\frac{\pi}{4} \left\{ \dfrac{1}{(2m-1)!} \left( \dfrac{\pi}{2} \right)^{2m-1} \right\}
				=(-1)^{m-1} \sum_{k=0}^{m-1} \left[ (-1)^{k} \lambda(2m-2k) \frac{1}{(2k)!} \left( \frac{\pi}{2} \right)^{2k} \right] \\
			&\frac{\pi}{4}  \left\{ \dfrac{1}{(2m)!} \left( \dfrac{\pi}{2} \right)^{2m} \right\}
				=(-1)^{m} \left[ \beta(2m+1)- \sum_{k=0}^{m-1} \left\{ (-1)^{k} \lambda(2m-2k) \frac{1}{(2k+1)!} \left( \frac{\pi}{2} \right)^{2k+1} \right\} \right]
		\end{align*}
		\end{proof}

\end{document}